\documentclass[11pt]{article}
\usepackage[margin=1in]{geometry}
\usepackage{amsmath}
\usepackage{amssymb}
\usepackage{amsthm}
\usepackage{graphicx}
\usepackage{booktabs}
\usepackage{tikz-cd}
\usepackage{multicol}
\usepackage{enumitem}
\usepackage{algorithm}
\usepackage{algpseudocode}
\algrenewcommand\algorithmicrequire{\textbf{Input:}}
\newtheorem{prop}{Proposition}
\newtheorem{cor}[prop]{Corollary}

\newtheorem{thm}[prop]{Theorem}

\usepackage[colorlinks=true,urlcolor=magenta]{hyperref}

\newcommand{\papertitle}{Multi-location Newsvendor: An Interpretable Inventory Policy}

\usepackage[authoryear,round]{natbib}

\begin{document}

\begingroup
\setlength{\parskip}{0pt}
\centering

{\LARGE\bfseries \papertitle\par}

\vspace{0.35cm}

{\large Soufiane Fafe\par}
\vspace{0.08cm}
{\small
Industrial Engineering and Operations Research,
Columbia University, New York, USA\par
\href{mailto:sf3159@columbia.edu}{sf3159@columbia.edu}\par}

\vspace{0.25cm}

{\large Vineet Goyal\par}
\vspace{0.08cm}
{\small
Industrial Engineering and Operations Research,
Columbia University, New York, USA\par
\href{mailto:vg2277@columbia.edu}{vg2277@columbia.edu}\par}

\endgroup

\vspace{0.3cm}

\begin{abstract}
We study a multilocation newsvendor problem in which inventory is chosen
before demand is observed and surplus stock can subsequently be transshipped
to locations with shortages. The objective is to minimize expected holding,
backlog, and transportation costs over a metric network. Under uniform
transportation costs, the optimality conditions balance local and aggregate
demand coverage and reduce to a scalar search under mild regularity
conditions. We use this structure to develop an interpretable iterative
heuristic for general metrics. The heuristic calibrates effective
transportation costs from lateral shipments under independent stocking and
uses them to construct location-dependent pooling weights. These weights
have a direct operational interpretation through the cost of serving demand
at other locations. Under the stated assumptions, the heuristic recovers
the exact uniform-metric solution. Numerical experiments on structured and
random networks show costs close to the sampled optimum on average and
substantial savings relative to independent stocking.
\end{abstract}

\section{Introduction}
\label{sec:introduction}

In many e-commerce systems, a seller distributes inventory across several warehouses to serve demand from different locations. Before demand is observed, the seller must decide how much inventory to place at each warehouse. After demand is realized, some locations may have excess inventory while others face shortages. Transshipping surplus inventory between locations allows the seller to meet demand that would otherwise remain unfulfilled.

Inventory decisions should therefore be planned jointly across the network. Treating each location as an independent newsvendor ignores the opportunity to pool inventory through transshipment. The value of this pooling depends on transportation costs: lower costs make sharing inventory across locations more economical, whereas higher costs increase the importance of protecting local demand. An inventory policy must account for this tradeoff together with the holding cost of unused stock and the backlog cost of unmet demand.

Sample-average approximation (SAA) is a classical approach to computing inventory allocations in this setting. It replaces expected cost with an average over sampled demand realizations and solves the resulting optimization problem. Our goal is to develop an interpretable stocking rule that makes explicit how transportation, holding, and backlog costs affect inventory decisions.

We begin with uniform transportation costs and characterize optimal inventory through a balance between local and aggregate demand coverage. Under mild regularity conditions, this characterization reduces computation to a scalar search. Furthermore, we build upon the structure of the optimal policy for the uniform costs to construct an iterative heuristic for general transportation metrics, with location-specific pooling parameters calibrated from lateral shipments. Numerical experiments compare the resulting policy with independent stocking and the SAA benchmark.

\subsection{Our Contributions}
\label{subsec:contributions}

Our main contributions are the following:

\medskip
\noindent\textbf{Optimal policy for a uniform metric.}
We first consider uniform transportation costs and characterize the optimal
inventory policy through a balance between local and aggregate demand
coverage. For a network of \(n\) locations, let \(F_j\) be the CDF of demand
at location \(j\), \(F_D\) the CDF of total network demand, and \(b\) and
\(h\) the unit backlog and holding costs. When the common lateral shipping
cost is \(\ell\), optimal inventories satisfy
\[
\alpha F_j(q_j^\star)
+
(1-\alpha)F_D\!\left(\sum_{i=1}^{n}q_i^\star\right)
=
\frac{b}{b+h},
\qquad
\alpha=\frac{\ell}{b+h}.
\]
The condition combines the probability that local inventory covers local
demand with the probability that total inventory covers total demand.
Transportation costs determine their relative weights. At \(\alpha=0\),
inventory can be shared freely and only the total quantity matters; at
\(\alpha=1\), the conditions reduce to independent newsvendor stocking.
For positive uniform costs, the optimal quantities can be computed through
a scalar search under the regularity conditions stated in
Section~\ref{sec:warmup-uniform}.

\medskip
\noindent\textbf{An iterative heuristic for general metrics.}
Building on this structure, we construct an interpretable policy that
replaces the common transportation cost and aggregate coverage term with
location-specific quantities. The resulting inventory vector, denoted
\(q^{\mathrm{ALG}}\), satisfies
\[
\alpha_j F_j(q_j^{\mathrm{ALG}})
+
(1-\alpha_j)
G_j\!\left(\sum_{i=1}^{n}w_{ij}q_i^{\mathrm{ALG}}\right)
=
\frac{b}{b+h}.
\]
Here, \(\alpha_j=\bar\ell_j/(b+h)\), where \(\bar\ell_j\) is the average
unit transportation cost from source \(j\). The CDF \(G_j\) describes weighted network demand, with
\(w_{ij}=\bar\ell_j/\ell_{ij}\) for \(i\ne j\) and \(w_{jj}=1\), where
\(\ell_{ij}\) is the unit shipping cost from \(j\) to \(i\). These weights
give greater emphasis to destinations that are cheaper to serve relative
to the calibrated reference cost. The policy therefore retains the balance
between local protection and pooling while accounting for transportation
geometry. Under positive uniform costs and well-defined calibration,
\(\bar\ell_j=\ell\) and all weights equal one, so the equations recover
the exact uniform-metric conditions.

\medskip
\noindent\textbf{Numerical performance and inventory placement.}
Our experiments show that the policy substantially reduces costs relative
to independent stocking while remaining close to the sampled optimum on
average. Across the reported metric classes, mean savings over independent
stocking range from \(19.64\%\) to \(53.34\%\), and mean gaps to SAA range
from \(2.24\%\) to \(10.06\%\). The savings arise from lower holding costs
that outweigh increased transportation and backlog costs. Comparisons of
inventory quantities and allocations show where the policy captures the
sampled optimal decisions and where differences remain. All comparisons
use a common demand sample, so the reported gaps are relative to the
sampled optimum.

\subsection{Related Work}
\label{subsec:related-work}

Our work relates to three streams of literature: inventory pooling
and network design, stocking with lateral transshipment, and
inventory placement and fulfillment.

\paragraph{Inventory pooling and network design.}
The newsvendor and inventory-pooling literature examines how
demand aggregation affects stocking decisions and expected costs
\citep{khouja1999single,eppen1979centralization,
corbett2006generalization}. The benefits depend on demand
variability and dependence, and greater variability need not
produce greater pooling savings
\citep{gerchak1992randomness,gerchak2003variability}.
Further analyses examine how these benefits change across demand
distributions, network sizes, and tail behavior, including settings
with uncertain yield
\citep{berman2011benefits,bimpikis2016pooling,mak2014pooling}.

Accounting for transportation costs connects pooling to the
configuration of the distribution network. Centralization and
location-inventory models weigh inventory savings against shipping
costs, with facility selection and customer assignments included
when the network is endogenous
\citep{chang1991centralization,daskin2002inventory,
shen2003joint}. Specialized optimization methods address the
resulting stochastic transportation-inventory design problem
\citep{shu2005stochastic}. These models determine where facilities
should operate and which demand locations they should serve.
We consider inventory decisions on a given network.

\paragraph{Stocking with lateral transshipment.}
With locations fixed, lateral transshipment links the inventories
held at different locations. Two-location and retail-network models
characterize stocking and transfer rules and their implications
for costs and service levels
\citep{das1975redistribution,tagaras1989pooling,wee2005optimal}.
Multiperiod and heterogeneous-location models extend these
decisions to repeated replenishment and larger systems, using
structural policies and network-flow-based computation
\citep{robinson1990optimal,herer2006multilocation}.
The distinction between proactive redistribution and reactive
transshipment organizes much of this literature; see
\citet{paterson2011review} for a review.

Resource-investment and newsvendor-network formulations provide
a broader framework for relating initial resource commitments
to operational decisions made after uncertainty is resolved
\citep{harrison1999multiresource,vanmieghem2002newsvendor}.
In particular, the newsvendor-network framework interprets
optimality conditions through multivariate critical fractiles.
For dynamic inventory systems, approximation methods and
simulation-based sensitivity estimates support the computation
of stocking policies
\citep{federgruen1984approximations,glasserman1995sensitivity}.

Calibration offers a related approach to constructing simpler
decision rules. \citet{archibald2010calibrations} uses
location-level calibrations to guide transshipment decisions
during a replenishment cycle. Our calibration uses lateral
shipments at a reference allocation to construct the initial
stocking equations. Decentralized and behavioral models introduce
additional considerations, including local stocking incentives,
inventory-sharing agreements, and reference-dependent ordering
\citep{rudi2001two,anupindi2001framework,ho2010reference}.

\paragraph{Inventory placement and fulfillment.}
Online retailing links inventory placement to a sequence of
fulfillment decisions. Linear-programming-based approaches account
for the future consequences of consuming inventory and the costs
of splitting multi-item orders
\citep{acimovic2015fulfillment,jasin2015rounding}.
Joint inventory and fulfillment formulations extend this analysis
to omnichannel demand and data-driven decisions about
inventory-holding locations and pooling arrangements
\citep{govindarajan2021omnichannel,wang2025placement}.
These models incorporate operating decisions throughout the
selling horizon, while our setting permits reactive transshipment
after demand is observed.

The information available about demand provides another
distinction between models. Moment-based approaches yield
distribution-free stocking rules, refinements that capture demand
asymmetry, and closed-form policies for multiperiod systems
\citep{gallego1993distribution,natarajan2018asymmetry,
mamani2017closedform}. For multiple locations, distributional
ambiguity has been addressed through pooling bounds and stocking
heuristics, as well as hierarchical policies with approximation
guarantees for sequential fulfillment on metrics
\citep{govindarajan2021pooling,foussoul2025metric}.
Our analysis assumes a specified joint demand distribution,
with the full demand vector observed before transshipment.

\vspace{2mm}
\noindent {\bf Outline.} Section~\ref{sec:model} defines the model and decomposes recourse cost into
aggregate mismatch and transportation costs. Section~\ref{sec:warmup-uniform}
develops the uniform-metric analysis. Section~\ref{sec:general-metric}
presents ALG and its uniform-metric exactness result.
Section~\ref{sec:numerics} reports the experiments, with additional results
and implementation details in Appendix~\ref{app:numerics}.
Section~\ref{sec:conclusion} concludes.

\section{Model and Preliminaries}\label{sec:model}

We consider a network of nodes $V=\{1,\dots,n\}$ with stochastic nonnegative demands
$(\xi_1,\dots,\xi_n)\sim F$, where $F$ is the joint cumulative distribution
function (CDF),
\[
\xi_i\ge0\ \text{a.s.},
\qquad
\mathbb E[\xi_i]<\infty
\quad\text{for every }i.
\]
The decision is an inventory vector $q=(q_1,\dots,q_n)\in\mathbb R_+^n$,
chosen before demand is realized. Let $h$ and $b$ denote the unit holding
and backlog costs, respectively. We assume throughout that
\[
h>0,\qquad b>0.
\]
After demand is realized, let $y_{ij}\ge0$ denote the amount shipped from
node $j$ to node $i$, and let $\ell_{ij}$ be the corresponding unit shipping
cost. The variable $y_{jj}$ represents local fulfillment, with $\ell_{jj}=0$.

We assume that $\ell=(\ell_{ij})$ is a metric on the set of inventory
locations $V$:\footnote{We allow zero transportation costs between distinct
locations, so the model also includes pseudometrics.}
\[
\ell_{ii}=0,\qquad
\ell_{ij}=\ell_{ji}\ge0,\qquad
\ell_{ij}\le\ell_{ik}+\ell_{kj}
\quad\text{for all }i,j,k\in V.
\]
Thus, transportation costs are symmetric, and direct shipment is no more
expensive than shipment through an intermediate location. We also assume that
\begin{equation}\label{eq:model-dominance}
\ell_{ij}\le h+b
\qquad \forall i\neq j.
\end{equation}
This assumption implies that, for any $i\neq j$, if there is unmet demand
at $i$ and a unit of excess inventory at $j$, then shipping the unit is
weakly no more expensive than paying holding at $j$ and backlog at $i$.

For any demand vector $\xi$, the recourse cost is
\begin{equation}
\label{eq:model-recourse}
\begin{aligned}
f(q,\xi)=\min_{y\ge 0}\Bigg\{&
\underbrace{h\sum_{j=1}^n\left(q_j-\sum_{i=1}^n y_{ij}\right)}_{\text{holding}}
+\underbrace{b\sum_{i=1}^n\left(\xi_i-\sum_{j=1}^n y_{ij}\right)}_{\text{backlog}}\\
&+\underbrace{\sum_{i\neq j}\ell_{ij}y_{ij}}_{\text{transport}}
:\ \sum_{i=1}^n y_{ij}\le q_j\ \forall j,
\qquad
\sum_{j=1}^n y_{ij}\le \xi_i\ \forall i
\Bigg\}.
\end{aligned}
\end{equation}
The constraints ensure that shipments do not exceed available inventory
at any source or realized demand at any destination.
The goal is to choose $q=(q_1,\dots,q_n)$ so as to minimize the expected total cost,
\begin{equation}\label{eq:model-value}
v(q):=\mathbb{E}[f(q,\xi)],
\qquad
\mathrm{OPT}:=\inf_{q\in\mathbb{R}_+^n} v(q).
\end{equation}

Let $Q:=\sum_{i=1}^n q_i$ and $D:=\sum_{i=1}^n \xi_i$ denote total
inventory and total demand. For convenience, define
\[
s_i:=(q_i-\xi_i)_+,\qquad
e_i:=(\xi_i-q_i)_+,\qquad
S:=\sum_{i=1}^n s_i,\qquad
E:=\sum_{i=1}^n e_i,
\]
where $(x)_+:=\max\{x,0\}$. Here, $s_i$ and $e_i$ are the surplus
inventory and unmet demand at node $i$ after local fulfillment, while
$S$ and $E$ are their network totals.

The metric assumption permits local fulfillment first without increasing
cost. Under \eqref{eq:model-dominance}, there exists an optimal recourse
solution that then ships $\min\{S,E\}$ units between distinct nodes.
In what follows, we select such an optimal solution when defining the
induced transportation cost.

We define the induced transportation cost by
\begin{equation}\label{eq:model-transport}
\begin{aligned}
T(q,\xi):=\min\Bigg\{&
\sum_{i\neq j}\ell_{ij}y_{ij}:
y\ge 0,
\quad
\sum_{i\neq j} y_{ij}\le s_j\ \forall j,
\quad
\sum_{j\neq i} y_{ij}\le e_i\ \forall i,
\\
&\hspace{38mm}
\sum_{i\neq j}y_{ij}=\min\{S,E\}
\Bigg\}.
\end{aligned}
\end{equation}
Therefore,
\begin{equation}\label{eq:model-decomp}
f(q,\xi)=h\big(Q-D\big)_+ + b\big(D-Q\big)_+ + T(q,\xi).
\end{equation}
The first two terms represent the remaining aggregate surplus or shortage
after transshipment. The final term accounts for the transportation cost
of reallocating inventory across locations.

\section{Warmup: Uniform Metric}\label{sec:warmup-uniform}

In this section, we consider the special case of a uniform metric.
In particular, $\ell_{ij}=\ell$ for all $i\neq j$. In this setting,
we show that the optimal inventory solution is characterized by a
fixed-point condition: a convex combination of the service level at
a node and the service level of cumulative inventory equals the
newsvendor fractile $b/(b+h)$. Moreover, the first-order conditions
provide useful structural insights into the optimal solution.
This case serves as a warm-up for Section~\ref{sec:general-metric}:
it provides an exact benchmark whose local--pooled structure
motivates the general-metric heuristic.

For the uniform metric, the shipping cost can be written in closed
form. In particular, the shipping cost depends only on the number
of shipped units, rather than on their sources and destinations.
Therefore, the recourse cost can be written as
\begin{equation}\label{eq:warmup-uniform-decomp}
f(q,\xi)
=
\ell\min\{S,E\}
+h(Q-D)_+
+b(D-Q)_+.
\end{equation}

The following theorem makes this local--pooled structure explicit:
optimal inventory balances the probability of covering demand at
each node with the probability of covering total network demand.

\begin{thm}[Optimal inventory under the uniform metric]\label{thm:warmup-uniform-opt}
Assume that each marginal CDF $F_i$ is continuous and that the total demand $D$ has continuous CDF $F_D$. Any optimum $q^\star$ satisfies
\begin{equation}\label{eq:warmup-uniform-kkt}
\frac{\ell}{h+b}\,F_j(q_j^\star)
+
\left(1-\frac{\ell}{h+b}\right)F_D(Q^\star)
=
\frac{b}{b+h},
\qquad j=1,\dots,n.
\end{equation}

\end{thm}

\begin{proof}
Taking expectation in \eqref{eq:warmup-uniform-decomp} gives
\[
v(q)
=
\ell\sum_{i=1}^n \mathbb{E}(\xi_i-q_i)_+
+h\mathbb{E}(Q-D)_+
+(b-\ell)\mathbb{E}(D-Q)_+.
\]
Differentiating with respect to $q_j$ yields
\[
\frac{\partial}{\partial q_j}\mathbb{E}(\xi_j-q_j)_+=F_j(q_j)-1,
\]
\[
\frac{\partial}{\partial q_j}\mathbb{E}(Q-D)_+=F_D(Q),
\qquad
\frac{\partial}{\partial q_j}\mathbb{E}(D-Q)_+=F_D(Q)-1.
\]
Combining these terms gives
\[
\frac{\partial v}{\partial q_j}(q)
=
\ell\,F_j(q_j)+(h+b-\ell)\,F_D(Q)-b.
\]
Moreover, the gradient map is monotone. For any \(q,r\in\mathbb R_+^n\), with
\(Q=\sum_j q_j\) and \(R=\sum_j r_j\),
\[
(q-r)^\top(\nabla v(q)-\nabla v(r))
=
\ell\sum_{j=1}^n (q_j-r_j)\big(F_j(q_j)-F_j(r_j)\big)
+
(h+b-\ell)(Q-R)\big(F_D(Q)-F_D(R)\big)
\ge0.
\]
Thus \(v\) is convex, since CDFs are nondecreasing and \(0\le\ell\le h+b\).

Define
\[
g_j(q):=\ell F_j(q_j)+(h+b-\ell)F_D(Q)-b.
\]
Since the feasible set is \(\mathbb R_+^n\), any optimum \(q^\star\) satisfies
\[
g_j(q^\star)\ge0,
\qquad
q_j^\star g_j(q^\star)=0,
\qquad j=1,\dots,n.
\]

If \(\ell=0\), then the objective depends on \(q\) only through \(Q\):
\[
v(q)=h\mathbb E(Q-D)_+ + b\mathbb E(D-Q)_+.
\]
Since \(F_D\) is continuous, any optimal total inventory \(Q^\star\) satisfies
\[
F_D(Q^\star)=\frac{b}{h+b},
\]
which gives \eqref{eq:warmup-uniform-kkt}. Hence assume \(\ell>0\).

Because demands are nonnegative and the CDFs are continuous, \(F_j(0)=0\) for all
\(j\), and \(F_D(0)=0\). Hence \(g_j(\mathbf 0)=-b<0\), so \(q^\star\neq\mathbf 0\).

Suppose, toward contradiction, that \(q_j^\star=0\) for some \(j\). There exists
\(k\) such that \(q_k^\star>0\) and \(F_k(q_k^\star)>0\). Otherwise,
\(F_i(q_i^\star)=0\) for every \(i\), which implies \(D>Q^\star\) a.s. and hence
\(F_D(Q^\star)=0\). But for any positive coordinate \(k\), complementarity gives
\(g_k(q^\star)=0\), while the previous identities would imply \(g_k(q^\star)=-b\),
a contradiction.

For such a coordinate \(k\), \(g_k(q^\star)=0\) gives
\[
(h+b-\ell)F_D(Q^\star)-b=-\ell F_k(q_k^\star)<0.
\]
At the zero coordinate \(j\),
\[
g_j(q^\star)
=
\ell F_j(0)+(h+b-\ell)F_D(Q^\star)-b
<0,
\]
contradicting \(g_j(q^\star)\ge0\). Therefore \(q_j^\star>0\) for every \(j\).
By complementarity, \(g_j(q^\star)=0\) for all \(j\). Dividing by \(h+b\) yields
\eqref{eq:warmup-uniform-kkt}.
\end{proof}

The following observations explain the structure of these optimality
conditions and their role in motivating the general-metric policy.
Let $\alpha=\ell/(b+h)$.
\begin{enumerate}[label=(\roman*),leftmargin=*]
    \item Equation \eqref{eq:warmup-uniform-kkt} can be written as
    \[
    \alpha F_j(q_j^\star)+(1-\alpha)F_D(Q^\star)=\frac{b}{b+h},
    \]
    so the uniform case is a convex combination of a local regime and a pooled regime.
    The term $F_j(q_j^\star)$ is the familiar node-wise service-level object from the
    classical newsvendor model: it asks whether node $j$ has enough stock for its own
    demand. By contrast, $F_D(Q^\star)$ is a system-wide service-level object: it asks
    whether the network as a whole has enough stock after demand is realized. The
    uniform metric is special because every node participates symmetrically in pooling,
    so the global contribution enters only through the scalar aggregate term
    $F_D(Q^\star)$.
    \item The two limit regimes clarify the meaning of $\alpha$. If $\alpha=1$
    (equivalently, $\ell=h+b$), then $F_j(q_j^\star)=\frac{b}{b+h}$ and the problem
    reduces to the independent policy. If $\alpha=0$ (equivalently, $\ell=0$), then
    $F_D(Q^\star)=\frac{b}{b+h}$ and only total inventory matters, since stock can be
    reallocated frictionlessly after demand. Thus $\alpha$ measures where the system lies
    between stand-alone operation and perfect pooling, a limit structure that the
    general-metric policy should also recover.
    \item With heterogeneous shipping costs, the pooled term should
    also depend on location. In other words, a unit stocked at node
    $j$ need not provide the same recourse value as a unit stocked at
    node $k$ when the transportation geometry is uneven. This
    observation suggests replacing the single pooled scalar from
    the uniform case with a location-dependent pooling term.
    Section~\ref{sec:general-metric} follows this construction,
    using effective transportation costs calibrated from lateral
    transshipments.
\end{enumerate}

The optimality conditions \eqref{eq:warmup-uniform-kkt} define a
coupled system for $q^\star$, but they can be solved by a
one-dimensional search under the following invertibility assumption.

\begin{cor}\label{cor:warmup-how-to-solve}
Under the assumptions of Theorem~\ref{thm:warmup-uniform-opt},
suppose $\ell>0$ and each marginal CDF $F_j$ has a continuous
inverse on $(0,1)$. Then the optimal inventory is given by
\[
q_j^\star=F_j^{-1}(r),
\qquad j=1,\dots,n,
\]
where $r\in(0,1)$ is the unique solution of
\[
\frac{\ell}{h+b}\,r
+
\left(1-\frac{\ell}{h+b}\right)
F_D\!\left(\sum_{j=1}^n F_j^{-1}(r)\right)
=
\frac{b}{b+h}.
\]
\end{cor}

\begin{proof}
Subtracting the $i$-th and $j$-th equations in
\eqref{eq:warmup-uniform-kkt} gives
\[
F_i(q_i^\star)=F_j(q_j^\star),
\qquad i,j=1,\dots,n.
\]
Let $r$ denote this common value. It lies in $(0,1)$:
if $r=0$, then $D>Q^\star$ almost surely, whereas if $r=1$,
then $D\le Q^\star$ almost surely. Either case contradicts
\eqref{eq:warmup-uniform-kkt}.
Invertibility therefore gives
\[
q_j^\star=F_j^{-1}(r),
\qquad j=1,\dots,n.
\]
Substituting this representation into
\eqref{eq:warmup-uniform-kkt} yields the stated scalar equation.

Its left-hand side is continuous and strictly increasing:
the quantiles and $F_D$ are nondecreasing, and $\ell>0$.
As $r$ approaches zero or one, the marginal quantiles approach
their lower or upper support endpoints, respectively. Continuity
of $F_D$ implies that the left-hand side approaches zero or one.
Thus the scalar equation has a unique solution, which can be
found by bisection. The resulting inventory vector has zero
gradient for the convex objective and is therefore optimal.
\end{proof}

\section{General-Metric Iterative Heuristic}
\label{sec:general-metric}

In this section, we propose an inventory policy for general transportation
metrics. The policy follows the local-versus-pooled structure of the
uniform-metric optimality conditions. We replace the common transportation
cost by a source-specific effective cost and the systemwide pooling term
by a weighted pooling term.

The effective costs are estimated from optimal lateral transshipments
under the independent newsvendor inventory allocation. They provide a
reference for measuring the cost of serving other locations from each
source. We first present the resulting iterative heuristic and explain
its construction. We then show that it recovers the exact optimality
conditions under a uniform metric.

\subsection{Construction of the Iterative Heuristic}
\label{subsec:iterative-heuristic}

The heuristic starts from the independent newsvendor quantities and uses
the resulting recourse flows to construct its pooling parameters. These
parameters remain fixed when solving the final inventory equations.

\begin{algorithm}[htbp]
\caption{General-metric iterative heuristic}
\label{alg:iterative-heuristic}
\begin{algorithmic}[1]
\Require Demand distribution \(F\), marginal CDFs \(F_j\), transportation
costs \(\ell_{ij}\), and costs \(b,h\).

\State \textbf{Initialize.} For each node \(j\), set the independent
newsvendor quantity
\begin{equation}
\label{eq:ih-initialization}
q_j^{(0)}
=
F_j^{-1}\left(\frac{b}{b+h}\right).
\end{equation}

\State \textbf{Evaluate recourse.} For each realization \(\xi\), solve the
transportation problem \eqref{eq:model-transport} at \(q^{(0)}\), and denote
the optimal lateral shipments by \(y_{ij}^{(0)}(\xi)\), \(i\neq j\).

\State \textbf{Estimate effective costs.} For every source node \(j\), compute
\begin{equation}
\label{eq:ih-effective-cost}
\bar{\ell}_j
=
\mathbb E\left[
\left.
\frac{
\sum_{i\neq j}y_{ij}^{(0)}(\xi)\ell_{ij}
}{
\sum_{i\neq j}y_{ij}^{(0)}(\xi)
}
\right|
\sum_{i\neq j}y_{ij}^{(0)}(\xi)>0
\right].
\end{equation}

\State \textbf{Construct pooling parameters.} For every \(i,j\), define
\begin{equation}
\label{eq:ih-alpha}
\alpha_j
=
\frac{\bar{\ell}_j}{b+h},
\end{equation}
\begin{equation}
\label{eq:ih-weights}
w_{ij}
=
\begin{cases}
\dfrac{\bar{\ell}_j}{\ell_{ij}}, & i\neq j,\\[8pt]
1, & i=j,
\end{cases}
\end{equation}
and
\begin{equation}
\label{eq:ih-weighted-demand}
Z_j
=
\sum_{i=1}^n w_{ij}\xi_i,
\end{equation}
with \(G_j\) denoting the CDF of \(Z_j\).

\State \textbf{Solve the fixed-point system.} Return a vector
\(q^{\mathrm{ALG}}\in\mathbb R_+^n\) satisfying
\begin{equation}
\label{eq:ih-system}
\alpha_jF_j(q_j)
+
(1-\alpha_j)
G_j\left(\sum_{i=1}^n w_{ij}q_i\right)
=
\frac{b}{b+h},
\qquad j=1,\dots,n.
\end{equation}
\end{algorithmic}
\end{algorithm}

The quantity \(\bar{\ell}_j\) represents the effective transportation cost
of using inventory at node \(j\) to satisfy demand elsewhere. For each
realization in which node \(j\) makes a positive lateral shipment, we
compute the transportation cost per unit shipped and then average over
these realizations. Thus \(\bar{\ell}_j\) accounts for the destinations
served under the independent allocation and the costs of reaching them.
Local fulfillment is excluded from this calculation.

The coefficient \(\alpha_j=\bar{\ell}_j/(b+h)\) compares this effective
transportation cost with the cost of holding a unit of surplus inventory
and leaving a unit of demand unmet. Since \(\ell_{ij}\le b+h\), we have
\(0\le\alpha_j\le1\). A smaller \(\alpha_j\) gives greater weight to pooling,
while a larger value gives greater weight to local demand protection.

For each destination \(i\neq j\), the weight
\(w_{ij}=\bar{\ell}_j/\ell_{ij}\) compares its transportation cost with the
effective cost from source \(j\). Destinations that can be served more
cheaply than this reference cost receive greater weight in the pooling
term, whereas more expensive destinations receive less weight. The
normalization \(w_{jj}=1\) leaves local inventory and demand unscaled.

Applying the same weights to inventory and demand gives the pooled
service term
\[
G_j\left(\sum_{i=1}^n w_{ij}q_i\right)
=
\mathbb P\left(
\sum_{i=1}^n w_{ij}\xi_i
\le
\sum_{i=1}^n w_{ij}q_i
\right).
\]
The heuristic balances this weighted pooling term with the local service
level \(F_j(q_j)\), following the structure of the uniform-metric
optimality condition.

\subsection{Exactness Under a Uniform Metric}
\label{subsec:ih-uniform-exactness}

When all lateral transportation costs are equal, the effective costs
coincide with the common transportation cost and every pooling weight
equals one. The following proposition shows that the heuristic therefore
recovers the uniform-metric solution.

\begin{prop}[Uniform-metric exactness]
\label{prop:ih-uniform-exactness}
Suppose that \(\ell_{ij}=\ell>0\) for every \(i\neq j\), that the marginal
CDFs and total-demand CDF are continuous, and that the conditional
expectations defining \(\bar{\ell}_j\) are well-defined. Then
\eqref{eq:ih-system} coincides with the uniform-metric optimality system
\eqref{eq:warmup-uniform-kkt}. Consequently, every solution of
\eqref{eq:ih-system} is optimal for the uniform-metric problem.
\end{prop}

\begin{proof}
Whenever source \(j\) makes a positive lateral shipment,
\[
\frac{
\sum_{i\neq j}y_{ij}^{(0)}(\xi)\ell_{ij}
}{
\sum_{i\neq j}y_{ij}^{(0)}(\xi)
}
=\ell.
\]
Hence \(\bar{\ell}_j=\ell\), which gives
\[
\alpha_j=\frac{\ell}{b+h},
\qquad
w_{ij}=1
\quad\text{for all }i,j.
\]
Therefore \(Z_j=D\), \(G_j=F_D\), and
\(\sum_i w_{ij}q_i=Q\). The heuristic equations become
\[
\frac{\ell}{b+h}F_j(q_j)
+
\left(1-\frac{\ell}{b+h}\right)F_D(Q)
=
\frac{b}{b+h},
\qquad j=1,\dots,n.
\]
These are the optimality conditions established in
Section~\ref{sec:warmup-uniform}. Since the uniform-metric objective is
convex, every feasible vector satisfying them is optimal.
\end{proof}

This result shows that the effective-cost calibration preserves the exact
uniform-metric benchmark. For general metrics, the next section evaluates
the heuristic's inventory allocations and costs against independent
stocking and the SAA benchmark.

\section{Numerical Experiments}
\label{sec:numerics}

We compare the general-metric iterative heuristic (ALG) with
independent stocking (IND) and sample-average approximation (SAA).
The experiments examine how much of the benefit from joint stocking
ALG captures and where gaps to the sampled optimum remain. We first
compare costs across metric classes, then examine transportation
scale and inventory placement to interpret the differences.

\subsection{Experimental setup}
\label{subsec:numerics-setup}

For each instance, ALG calibration, its pooling CDFs, SAA, and policy
evaluation use the same \(N=2000\) demand realizations. The sampled
cost is
\[
\widehat v(q)=\frac1N\sum_{s=1}^{N}f(q,\xi^{(s)}).
\]
IND selects the analytic marginal newsvendor quantities
\citep{arrow1951optimal},
\[
q_j^{\mathrm{IND}}
=F_j^{-1}\left(\frac{b}{b+h}\right).
\]
ALG follows Algorithm~\ref{alg:iterative-heuristic}.
SAA solves the sampled optimization problem
\citep{shapiro2009lectures},
\[
q^{\mathrm{SAA}}\in\arg\min_{q\ge0}\widehat v(q),
\]
and provides the sampled optimum to the numerical tolerances
specified in Appendix~\ref{app:numerics}. All policies use optimal
transshipment recourse during evaluation.

We report the sampled optimality gap and improvement over IND:
\[
\Delta(P)
=100\frac{\widehat v(q^P)-\widehat v(q^{\mathrm{SAA}})}
{\widehat v(q^{\mathrm{SAA}})},
\qquad
\mathrm{Imp}(P;\mathrm{IND})
=100\frac{\widehat v(q^{\mathrm{IND}})-\widehat v(q^P)}
{\widehat v(q^{\mathrm{IND}})}.
\]
We also report total inventory relative to SAA,
\(Q/Q^{\mathrm{SAA}}\). Percentages and inventory ratios are
calculated per instance before averaging; P90 denotes the
90th percentile of instance gaps. Holding, backlog, and
transportation costs are the sample averages of the three terms
in \eqref{eq:model-decomp}.

All summaries use instances with a certified SAA solution and a
verified nonnegative ALG solution. Within each comparison, the
policies are evaluated on the same instances. The reported gaps
therefore measure performance relative to the optimum of the
shared sampled problem.

\begin{table}[htbp]
\centering
\small
\caption{Experimental parameter grid.}
\label{tab:numerics-parameters}
\begin{tabular}{ll}
\toprule
Parameter & Values \\
\midrule
Number of nodes \(n\) & \(10,20,30,40,50\) \\
Mean demand & \(\mu=1\) \\
Demand laws & Truncated Gaussian, log-normal, gamma, Pareto \\
Coefficient of variation & \(\mathrm{CV}\in\{0.25,0.50\}\) \\
Critical fractile & \(b/(b+h)\in\{0.80,0.85,0.90,0.95\}\) \\
Combined mismatch cost & \(b+h=5\) \\
Transportation scale \(\gamma\) & \(0.5,1,2,3,4.5\) \\
Sample size \(N\) & \(2000\) \\
\bottomrule
\end{tabular}
\end{table}

Table~\ref{tab:numerics-parameters} gives the parameter grid.
Setting \(b+h=5\) normalizes costs: multiplying holding, backlog,
and transportation costs by the same positive factor leaves
optimal quantities and ALG pooling parameters unchanged.
The relevant quantities are \(b/(b+h)\) and transportation costs
relative to \(b+h\). The selected fractiles examine increasing
shortage penalties, with \(b/h\) ranging from \(4\) to \(19\).

Demands are independent and identically distributed across nodes.
The four laws have matching means and CVs, allowing us to compare
distributional shape while controlling the first two moments.
Gaussian demand is truncated below at zero, with moments matched
after truncation; Pareto demand follows the Type-I distribution.
The two CVs represent standard deviations equal to one-quarter
and one-half of mean demand. Setting \(\mu=1\) normalizes demand
units. Distribution parameters appear in
Table~\ref{tab:app-demand-parameters}.

For base distances \(d_{ij}\), transportation costs are
\[
\ell_{ij}^{(\gamma)}
=\gamma\frac{d_{ij}}{\max_{k\ne l}d_{kl}},
\qquad i\ne j,
\qquad \ell_{ii}^{(\gamma)}=0.
\]
The maximum base distance is positive. This normalization preserves
relative distances while setting the largest transportation cost
between \(10\%\) and \(90\%\) of \(b+h\). The grid therefore ranges
from inexpensive transshipment to costs close to the combined
holding and backlog cost.

Structured instances are stars with hub--leaf distance \(1/2\)
and leaf--leaf distance \(1\), lines with \(d_{ij}=|i-j|\), and
two equal clusters with within-group distance \(0.1\) and
between-group distance \(1\). These geometries represent a common
hub, a spatial ordering, and inexpensive sharing within groups.
Random instances use Euclidean distances between independent
uniform locations in \([0,1]^2\), or unit-edge shortest-path
distances in \(G(n,0.4)\) conditioned on connectivity.

Appendix~\ref{app:numerics} reports implementation details,
uniform-metric validation, cost decompositions, additional
sensitivity analyses, and construction times.

\subsection{Overall performance}
\label{subsec:numerics-performance}

ALG captures substantial cost savings from joint stocking.
Across the five metric classes, mean improvements over IND range
from \(19.64\%\) to \(53.34\%\), while mean gaps to SAA range from
\(2.24\%\) to \(10.06\%\). Table~\ref{tab:numerics-main} reports
the cost comparisons and total-inventory ratios.

\begin{table}[htbp]
\centering
\small
\caption{Policy performance on paired instances, \(N=2000\).
Gaps are relative to the sampled SAA optimum; improvements are
relative to IND.}
\label{tab:numerics-main}
\resizebox{\textwidth}{!}{%
\begin{tabular}{llrrrrrr}
\toprule
Metric & Policy & \(N\) & Paired cases
& Mean gap (\%) & P90 gap (\%)
& Imp. vs. IND (\%) & \(Q/Q^{\mathrm{SAA}}\) \\
\midrule
Star & SAA & 2000 & 3998 & 0.00 & 0.00 & 26.09 & 1.000 \\
Star & IND & 2000 & 3998 & 44.36 & 113.81 & 0.00 & 1.179 \\
Star & ALG & 2000 & 3998 & 8.30 & 15.47 & 19.64 & 1.000 \\
\midrule
Line & SAA & 2000 & 3486 & 0.00 & 0.00 & 56.64 & 1.000 \\
Line & IND & 2000 & 3486 & 166.30 & 298.80 & 0.00 & 1.277 \\
Line & ALG & 2000 & 3486 & 10.06 & 23.46 & 53.34 & 1.030 \\
\midrule
Clustered & SAA & 2000 & 3573 & 0.00 & 0.00 & 53.03 & 1.000 \\
Clustered & IND & 2000 & 3573 & 139.45 & 258.25 & 0.00 & 1.285 \\
Clustered & ALG & 2000 & 3573 & 2.53 & 4.60 & 51.66 & 1.010 \\
\midrule
Geometric & SAA & 2000 & 7369 & 0.00 & 0.00 & 43.37 & 1.000 \\
Geometric & IND & 2000 & 7369 & 95.48 & 188.66 & 0.00 & 1.221 \\
Geometric & ALG & 2000 & 7369 & 6.33 & 15.68 & 40.32 & 1.009 \\
\midrule
Random graph & SAA & 2000 & 11816 & 0.00 & 0.00 & 33.78 & 1.000 \\
Random graph & IND & 2000 & 11816 & 64.67 & 151.11 & 0.00 & 1.243 \\
Random graph & ALG & 2000 & 11816 & 2.24 & 4.20 & 32.68 & 1.003 \\
\bottomrule
\end{tabular}%
}
\end{table}

Random graphs and clustered instances give the smallest mean
ALG gaps. Lines illustrate that large savings over IND can coexist
with a noticeable gap to SAA: ALG reduces cost by \(53.34\%\)
relative to IND but retains a mean sampled gap of \(10.06\%\).
Its P90 gap is also largest on lines, at \(23.46\%\).

The cost decomposition explains the savings over independent
stocking. In every metric class, ALG's lower mean holding cost
more than offsets its higher mean backlog and transportation
costs; see Table~\ref{tab:app-cost-components}. Inventory sharing
allows the policy to carry less surplus stock while accepting
additional shipments and shortages.

The inventory ratios support this interpretation. Mean total
inventory relative to SAA lies between \(1.000\) and \(1.030\)
for ALG, compared with \(1.179\) to \(1.285\) for IND. ALG therefore
produces aggregate quantities much closer to the sampled optimum.
The remaining cost gaps also depend on where this inventory is
placed, which we examine next.

\subsection{Transportation costs and inventory placement}
\label{subsec:numerics-interpretation}

The relationship between transportation scale and ALG accuracy
differs across metric classes. Table~\ref{tab:numerics-scale}
varies the maximum transportation cost while preserving the
relative distances of each base metric.

\begin{table}[htbp]
\centering
\small
\caption{Mean ALG sampled gap (\%) by transportation scale,
\(N=2000\), on paired instances at each scale.}
\label{tab:numerics-scale}
\resizebox{\textwidth}{!}{%
\begin{tabular}{rrrrrrr}
\toprule
\(\gamma\) & \(N\)
& Star & Line & Clustered & Geometric & Random graph \\
\midrule
0.50 & 2000 & 3.47 & 12.45 & 3.65 & 5.20 & 4.42 \\
1.00 & 2000 & 6.56 & 11.30 & 3.13 & 7.55 & 5.48 \\
2.00 & 2000 & 9.66 & 9.79 & 2.90 & 7.36 & 2.73 \\
3.00 & 2000 & 10.83 & 8.34 & 1.82 & 6.48 & 0.64 \\
4.50 & 2000 & 10.96 & 7.59 & 1.27 & 4.88 & 0.13 \\
\bottomrule
\end{tabular}%
}
\end{table}

On stars, the mean gap increases from \(3.47\%\) to \(10.96\%\)
over the tested scales. On lines and clustered instances, it
decreases from \(12.45\%\) to \(7.59\%\) and from \(3.65\%\) to
\(1.27\%\), respectively. Geometric and random-graph results are
nonmonotone. Increasing transportation costs therefore does not
produce a common direction of change in the reported gaps.
Because the paired subsets can differ across scales, these
patterns describe comparisons of the reported group means.

\begin{figure}[htbp]
\centering
\IfFileExists{inventory_allocation.png}{%
\includegraphics[width=\textwidth]{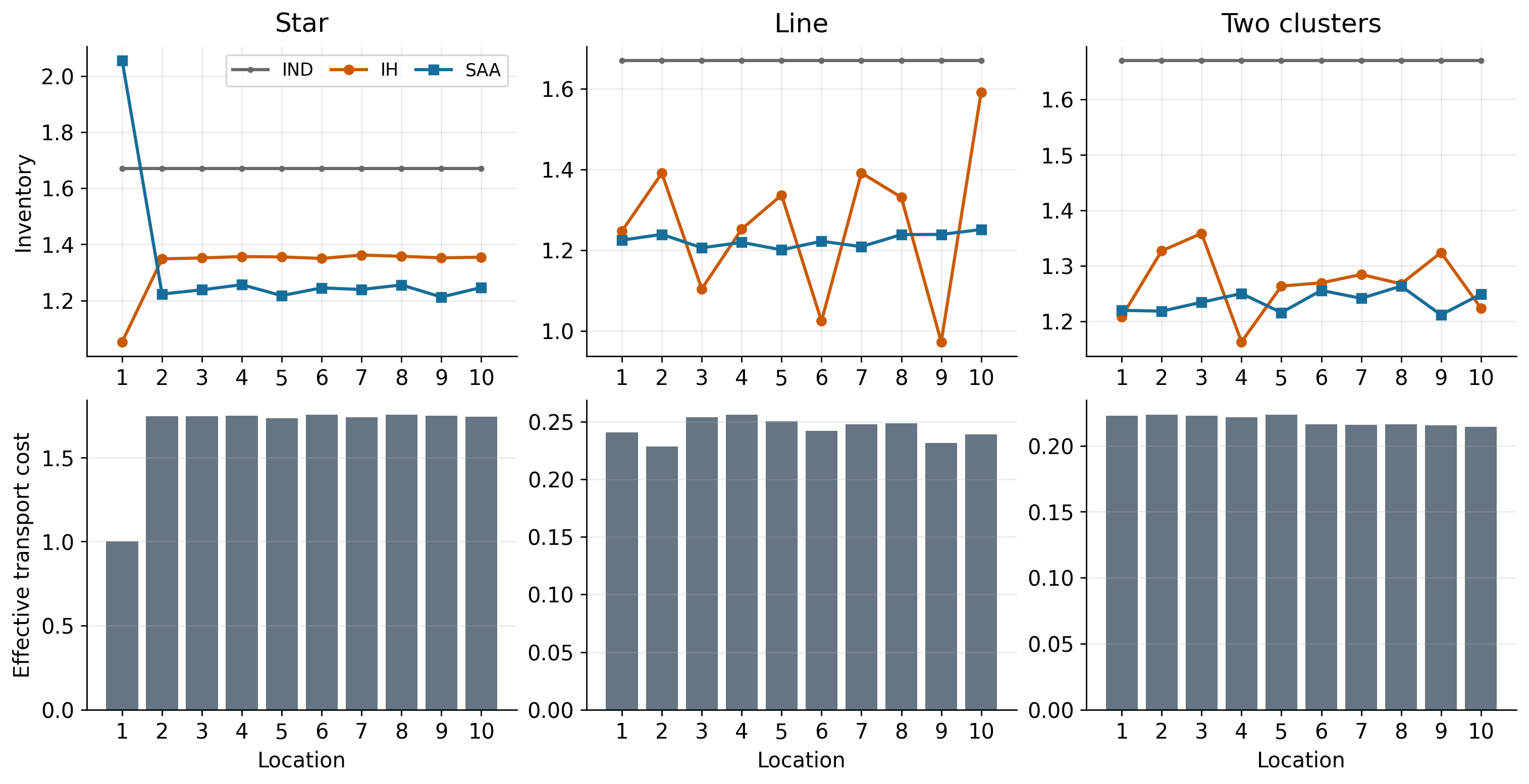}%
}{%
\fbox{\parbox{0.92\textwidth}{\centering
\textit{[Figure pending:
\texttt{inventory\_allocation.png}.]}}}%
}
\caption{Inventory allocations (top) and effective transportation
costs calibrated at IND (bottom) for star, line, and two-cluster
instances. Gamma demand, \(n=10\), \(N=2000\),
\(\mathrm{CV}=0.50\), \(b/(b+h)=0.90\), \(\gamma=2\).}
\label{fig:numerics-allocation}
\end{figure}

Inventory placement provides a second perspective on the gaps.
Figure~\ref{fig:numerics-allocation} compares node-level quantities
on three structured instances with common demand and cost
parameters. The lower panels show the effective transportation
costs calibrated at IND.

IND assigns the same quantity to every node because the marginal
demands are identical. ALG and SAA account for the network but
can produce different allocation profiles. In the line and
clustered examples, ALG exhibits greater node-to-node variation
than SAA. Thus, a total-inventory ratio close to one does not
ensure that the policy reproduces the sampled optimal placement.

This distinction matters because holding and backlog costs depend
on total inventory, whereas transportation costs also depend on
its allocation. For stars, transportation accounts for the largest
component of ALG's mean cost difference from SAA in
Table~\ref{tab:app-cost-components}. The aggregate inventory
comparisons and the allocation figure therefore assess different
aspects of the stocking decision.

\section{Conclusion}
\label{sec:conclusion}

We develop an interpretable stocking policy for the multilocation
newsvendor problem with reactive transshipment. Under uniform
transportation costs, the optimality conditions balance local
and aggregate demand coverage and reduce to a scalar search
under the stated regularity conditions. For general metrics,
ALG uses lateral shipments under independent stocking to
calibrate effective transportation costs and pooling weights.
With well-defined calibration, the construction recovers the
optimality conditions for positive uniform costs.

On the reported structured and random instances, ALG substantially
reduces mean costs relative to independent stocking and remains
close to the sampled optimum on average. Lower holding costs
more than offset higher backlog and transportation costs.
The smallest mean gaps occur on clustered instances and random
graphs. Lines produce both the largest savings over independent
stocking and the largest remaining mean gaps to SAA.
The allocation comparisons also show that total inventory close
to the sampled optimum can coexist with differences in placement,
which affect transportation costs.

Further work could identify conditions on the metric and demand
distribution under which ALG admits performance guarantees.
Evaluating the calibration under dependent or heterogeneous
demands would also extend the numerical analysis.

\section*{Declaration of generative AI and AI-assisted technologies}

During the preparation of this manuscript, the authors
used a large language model (LLM) to assist with grammar,
language clarity, LaTeX formatting, and manuscript
organization. The authors reviewed the output and take
full responsibility for the final content of the manuscript.

\bibliographystyle{plainnat}
\bibliography{references}

@article{corbett2006generalization,
  author  = {Corbett, Charles J. and Rajaram, Kumar},
  title   = {A Generalization of the Inventory Pooling Effect to Nonnormal Dependent Demand},
  journal = {Manufacturing \& Service Operations Management},
  volume  = {8},
  number  = {4},
  pages   = {351--358},
  year    = {2006},
  doi     = {10.1287/msom.1060.0117}
}

@article{bimpikis2016pooling,
  author  = {Bimpikis, Kostas and Markakis, Mihalis G.},
  title   = {Inventory Pooling Under Heavy-Tailed Demand},
  journal = {Management Science},
  volume  = {62},
  number  = {6},
  pages   = {1800--1813},
  year    = {2016},
  doi     = {10.1287/mnsc.2015.2204}
}

@article{herer2006multilocation,
  author  = {Herer, Yale T. and Tzur, Michal and Y{\"u}cesan, Enver},
  title   = {The Multilocation Transshipment Problem},
  journal = {IIE Transactions},
  volume  = {38},
  number  = {3},
  pages   = {185--200},
  year    = {2006},
  doi     = {10.1080/07408170500434539}
}

@article{archibald2010calibrations,
  author  = {Archibald, T. W. and Black, D. P. and Glazebrook, K. D.},
  title   = {The Use of Simple Calibrations of Individual Locations in Making Transshipment Decisions in a Multi-Location Inventory Network},
  journal = {Journal of the Operational Research Society},
  volume  = {61},
  number  = {2},
  pages   = {294--305},
  year    = {2010},
  doi     = {10.1057/jors.2008.127}
}

@article{anupindi2001framework,
  author  = {Anupindi, Ravi and Bassok, Yehuda and Zemel, Eitan},
  title   = {A General Framework for the Study of Decentralized Distribution Systems},
  journal = {Manufacturing \& Service Operations Management},
  volume  = {3},
  number  = {4},
  pages   = {349--368},
  year    = {2001},
  doi     = {10.1287/msom.3.4.349.9973}
}

@article{govindarajan2021pooling,
  author  = {Govindarajan, Aravind and Sinha, Amitabh and Uichanco, Joline},
  title   = {Distribution-Free Inventory Risk Pooling in a Multilocation Newsvendor},
  journal = {Management Science},
  volume  = {67},
  number  = {4},
  pages   = {2272--2291},
  year    = {2021},
  doi     = {10.1287/mnsc.2020.3719}
}

@misc{foussoul2025metric,
  author = {Foussoul, Ayoub and Goyal, Vineet},
  title  = {Distributionally Robust Newsvendor on a Metric},
  year   = {2025},
  note   = {arXiv:2410.12134v2, revised June 2, 2025},
  doi    = {10.48550/arXiv.2410.12134}
}

@article{arrow1951optimal,
  author  = {Arrow, Kenneth J. and Harris, Theodore and Marschak, Jacob},
  title   = {Optimal Inventory Policy},
  journal = {Econometrica},
  volume  = {19},
  number  = {3},
  pages   = {250--272},
  year    = {1951},
  doi     = {10.2307/1906813}
}

@article{eppen1979centralization,
  author  = {Eppen, Gary D.},
  title   = {Note---Effects of Centralization on Expected Costs in a Multi-Location Newsboy Problem},
  journal = {Management Science},
  volume  = {25},
  number  = {5},
  pages   = {498--501},
  year    = {1979},
  doi     = {10.1287/mnsc.25.5.498}
}

@article{das1975redistribution,
  author  = {Das, Chandrasekhar},
  title   = {Supply and Redistribution Rules for Two-Location Inventory Systems: One-Period Analysis},
  journal = {Management Science},
  volume  = {21},
  number  = {7},
  pages   = {765--776},
  year    = {1975},
  doi     = {10.1287/mnsc.21.7.765}
}

@article{robinson1990optimal,
  author  = {Robinson, Lawrence W.},
  title   = {Optimal and Approximate Policies in Multiperiod, Multilocation Inventory Models with Transshipments},
  journal = {Operations Research},
  volume  = {38},
  number  = {2},
  pages   = {278--295},
  year    = {1990},
  doi     = {10.1287/opre.38.2.278}
}

@article{rudi2001two,
  author  = {Rudi, Nils and Kapur, Sandeep and Pyke, David F.},
  title   = {A Two-Location Inventory Model with Transshipment and Local Decision Making},
  journal = {Management Science},
  volume  = {47},
  number  = {12},
  pages   = {1668--1680},
  year    = {2001},
  doi     = {10.1287/mnsc.47.12.1668.10235}
}

@article{paterson2011review,
  author  = {Paterson, Colin and Kiesm{\"u}ller, Gudrun and Teunter, Ruud and Glazebrook, Kevin},
  title   = {Inventory Models with Lateral Transshipments: A Review},
  journal = {European Journal of Operational Research},
  volume  = {210},
  number  = {2},
  pages   = {125--136},
  year    = {2011}
}

@article{vanmieghem2002newsvendor,
  author  = {Van Mieghem, Jan A. and Rudi, Nils},
  title   = {Newsvendor Networks: Inventory Management and Capacity Investment with Discretionary Activities},
  journal = {Manufacturing \& Service Operations Management},
  volume  = {4},
  number  = {4},
  pages   = {313--335},
  year    = {2002},
  doi     = {10.1287/msom.4.4.313.5728}
}

@article{wee2005optimal,
  author  = {Wee, Kwan Eng and Dada, Maqbool},
  title   = {Optimal Policies for Transshipping Inventory in a Retail Network},
  journal = {Management Science},
  volume  = {51},
  number  = {10},
  pages   = {1519--1533},
  year    = {2005},
  doi     = {10.1287/mnsc.1050.0441}
}

@book{shapiro2009lectures,
  author    = {Shapiro, Alexander and Dentcheva, Darinka and Ruszczynski, Andrzej},
  title     = {Lectures on Stochastic Programming: Modeling and Theory},
  publisher = {Society for Industrial and Applied Mathematics},
  address   = {Philadelphia},
  year      = {2009}
}

@article{khouja1999single,
  author  = {Khouja, M.},
  title   = {The Single-Period (News-Vendor) Problem: Literature Review and Suggestions for Future Research},
  journal = {Omega},
  volume  = {27},
  number  = {5},
  pages   = {537--553},
  year    = {1999},
  doi     = {10.1016/S0305-0483(99)00017-1}
}

@article{chang1991centralization,
  author  = {Chang, P.-L. and Lin, C.-T.},
  title   = {On the Effect of Centralization on Expected Costs in a Multi-Location Newsboy Problem},
  journal = {Journal of the Operational Research Society},
  volume  = {42},
  number  = {11},
  pages   = {1025--1030},
  year    = {1991},
  doi     = {10.1057/jors.1991.193}
}

@article{tagaras1989pooling,
  author  = {Tagaras, G.},
  title   = {Effects of Pooling on the Optimization and Service Levels of Two-Location Inventory Systems},
  journal = {IIE Transactions},
  volume  = {21},
  number  = {3},
  pages   = {250--257},
  year    = {1989},
  doi     = {10.1080/07408178908966229}
}

@article{gerchak1992randomness,
  author  = {Gerchak, Y. and Mossman, D.},
  title   = {On the Effect of Demand Randomness on Inventories and Costs},
  journal = {Operations Research},
  volume  = {40},
  number  = {4},
  pages   = {804--807},
  year    = {1992},
  doi     = {10.1287/opre.40.4.804}
}

@article{gerchak2003variability,
  author  = {Gerchak, Y. and He, Q.-M.},
  title   = {On the Relation Between the Benefits of Risk Pooling and the Variability of Demand},
  journal = {IIE Transactions},
  volume  = {35},
  number  = {11},
  pages   = {1027--1031},
  year    = {2003},
  doi     = {10.1080/07408170304399}
}

@article{berman2011benefits,
  author  = {Berman, O. and Krass, D. and Tajbakhsh, M. M.},
  title   = {On the Benefits of Risk Pooling in Inventory Management},
  journal = {Production and Operations Management},
  volume  = {20},
  number  = {1},
  pages   = {57--71},
  year    = {2011},
  doi     = {10.1111/j.1937-5956.2010.01134.x}
}

@article{mak2014pooling,
  author  = {Mak, H.-Y. and Shen, Z.-J. M.},
  title   = {Pooling and Dependence of Demand and Yield in Multiple-Location Inventory Systems},
  journal = {Manufacturing \& Service Operations Management},
  volume  = {16},
  number  = {2},
  pages   = {263--269},
  year    = {2014},
  doi     = {10.1287/msom.2013.0469}
}

@article{daskin2002inventory,
  author  = {Daskin, M. S. and Coullard, C. R. and Shen, Z.-J. M.},
  title   = {An Inventory-Location Model: Formulation, Solution Algorithm and Computational Results},
  journal = {Annals of Operations Research},
  volume  = {110},
  pages   = {83--106},
  year    = {2002},
  doi     = {10.1023/A:1020763400324}
}

@article{shen2003joint,
  author  = {Shen, Z.-J. M. and Coullard, C. R. and Daskin, M. S.},
  title   = {A Joint Location-Inventory Model},
  journal = {Transportation Science},
  volume  = {37},
  number  = {1},
  pages   = {40--55},
  year    = {2003},
  doi     = {10.1287/trsc.37.1.40.12823}
}

@article{shu2005stochastic,
  author  = {Shu, J. and Teo, C.-P. and Shen, Z.-J. M.},
  title   = {Stochastic Transportation-Inventory Network Design Problem},
  journal = {Operations Research},
  volume  = {53},
  number  = {1},
  pages   = {48--60},
  year    = {2005},
  doi     = {10.1287/opre.1040.0140}
}

@article{harrison1999multiresource,
  author  = {Harrison, J. M. and Van Mieghem, J. A.},
  title   = {Multi-Resource Investment Strategies: Operational Hedging Under Demand Uncertainty},
  journal = {European Journal of Operational Research},
  volume  = {113},
  number  = {1},
  pages   = {17--29},
  year    = {1999},
  doi     = {10.1016/S0377-2217(97)00425-6}
}

@article{federgruen1984approximations,
  author  = {Federgruen, A. and Zipkin, P.},
  title   = {Approximations of Dynamic, Multilocation Production and Inventory Problems},
  journal = {Management Science},
  volume  = {30},
  number  = {1},
  pages   = {69--84},
  year    = {1984},
  doi     = {10.1287/mnsc.30.1.69}
}

@article{glasserman1995sensitivity,
  author  = {Glasserman, P. and Tayur, S.},
  title   = {Sensitivity Analysis for Base-Stock Levels in Multiechelon Production-Inventory Systems},
  journal = {Management Science},
  volume  = {41},
  number  = {2},
  pages   = {263--281},
  year    = {1995},
  doi     = {10.1287/mnsc.41.2.263}
}

@article{ho2010reference,
  author  = {Ho, T.-H. and Lim, N. and Cui, T. H.},
  title   = {Reference Dependence in Multilocation Newsvendor Models: A Structural Analysis},
  journal = {Management Science},
  volume  = {56},
  number  = {11},
  pages   = {1891--1910},
  year    = {2010},
  doi     = {10.1287/mnsc.1100.1225}
}

@article{gallego1993distribution,
  author  = {Gallego, G. and Moon, I.},
  title   = {The Distribution Free Newsboy Problem: Review and Extensions},
  journal = {Journal of the Operational Research Society},
  volume  = {44},
  number  = {8},
  pages   = {825--834},
  year    = {1993},
  doi     = {10.1057/jors.1993.141}
}

@article{natarajan2018asymmetry,
  author  = {Natarajan, K. and Sim, M. and Uichanco, J.},
  title   = {Asymmetry and Ambiguity in Newsvendor Models},
  journal = {Management Science},
  volume  = {64},
  number  = {7},
  pages   = {3146--3167},
  year    = {2018},
  doi     = {10.1287/mnsc.2017.2773}
}

@article{mamani2017closedform,
  author  = {Mamani, H. and Nassiri, S. and Wagner, M. R.},
  title   = {Closed-Form Solutions for Robust Inventory Management},
  journal = {Management Science},
  volume  = {63},
  number  = {5},
  pages   = {1625--1643},
  year    = {2017},
  doi     = {10.1287/mnsc.2015.2391}
}

@article{acimovic2015fulfillment,
  author  = {Acimovic, J. and Graves, S. C.},
  title   = {Making Better Fulfillment Decisions on the Fly in an Online Retail Environment},
  journal = {Manufacturing \& Service Operations Management},
  volume  = {17},
  number  = {1},
  pages   = {34--51},
  year    = {2015},
  doi     = {10.1287/msom.2014.0505}
}

@article{jasin2015rounding,
  author  = {Jasin, S. and Sinha, A.},
  title   = {An {LP}-Based Correlated Rounding Scheme for Multi-Item Ecommerce Order Fulfillment},
  journal = {Operations Research},
  volume  = {63},
  number  = {6},
  pages   = {1336--1351},
  year    = {2015},
  doi     = {10.1287/opre.2015.1441}
}

@article{govindarajan2021omnichannel,
  author  = {Govindarajan, A. and Sinha, A. and Uichanco, J.},
  title   = {Joint Inventory and Fulfillment Decisions for Omnichannel Retail Networks},
  journal = {Naval Research Logistics},
  volume  = {68},
  number  = {6},
  pages   = {779--794},
  year    = {2021},
  doi     = {10.1002/nav.21969}
}

@article{wang2025placement,
  author  = {Wang, Y. and Minner, S.},
  title   = {Data-Driven Multi-Location Inventory Placement in Digital Commerce},
  journal = {Computers \& Industrial Engineering},
  volume  = {200},
  pages   = {110842},
  year    = {2025},
  doi     = {10.1016/j.cie.2024.110842}
}

\clearpage
\appendix
\section{Additional Numerical Results}
\label{app:numerics}

This appendix reports implementation details, uniform-metric validation,
cost decompositions, sensitivity analyses, and construction times.
Experiments use \(N=2000\) and the parameter grid in
Table~\ref{tab:numerics-parameters}. Unless varied explicitly, the
within-cluster base distance is \(0.1\) and the random-graph edge
probability is \(p=0.4\). Policy comparisons use the completed paired
instances defined in Section~\ref{sec:numerics}. Percentages are computed
per instance before averaging; the included instances can differ
between sensitivity settings.

\subsection{Implementation and demand parameters}
\label{app:numerics-implementation}

Structured settings use five sample replications. Random settings use
twenty network and sample realizations, reused across demand and cost
parameters.

Table~\ref{tab:app-demand-parameters} specifies the four demand laws,
each with mean \(\mu=1\) and the indicated CV. Gaussian moments are
matched after truncation below at zero. Pareto demand follows the
Type-I distribution.

\begin{table}[htbp]
\centering
\small
\caption{Demand parameters, rounded for presentation.}
\label{tab:app-demand-parameters}
\resizebox{\textwidth}{!}{%
\begin{tabular}{llll}
\toprule
Demand law & Parameter convention
& \(\mathrm{CV}=0.25\) & \(\mathrm{CV}=0.50\) \\
\midrule
Truncated Gaussian
& Underlying mean and standard deviation
& \((0.999966,0.250067)\) & \((0.948636,0.548966)\) \\
Log-normal
& Mean and standard deviation of log demand
& \((-0.030312,0.246221)\) & \((-0.111572,0.472381)\) \\
Gamma
& Shape and scale
& \((16,0.0625)\) & \((4,0.25)\) \\
Pareto
& Shape and minimum demand
& \((5.123106,0.804806)\) & \((3.236068,0.690983)\) \\
\bottomrule
\end{tabular}%
}
\end{table}

Effective transportation costs are calibrated at
\(q^{\mathrm{IND}}\) using only lateral shipments, as in
\eqref{eq:ih-effective-cost}. For each source, we average the
scenario-wise cost per unit shipped over realizations with export
volume exceeding \(10^{-10}\). Uniform, star, and clustered instances
use proportional allocations among tied optimal flows. Line and
random-graph calibrations average eight seeded optimal flow selections;
geometric instances use one optimal flow selection. The calibrated
parameters remain fixed during the ALG solve.

Marginal CDFs are analytic. Weighted-demand CDFs use linear interpolation
of the sorted common-sample values at midpoint probability levels
\((s-1/2)/N\), \(s=1,\ldots,N\). Uniform ALG uses a scalar bracketed
search. For general metrics, inverse-CDF and forward-CDF formulations
are solved numerically, with bounded residual least squares used when
needed. Reported ALG solutions are nonnegative, have maximum residual
at most \(10^{-6}\) in \eqref{eq:ih-system}, and have raw empirical-CDF
residual at most \(1/(2N)+10^{-6}\). The reported paired comparisons
use positive off-diagonal costs and the original inverse weights.

SAA uses an exact empirical marginal-cost procedure for uniform and
two-cluster instances, and LP or decomposition methods for the remaining
metrics. Every reported SAA solution is evaluated on all 2,000 scenarios.
For decomposition, the feasible cost and a valid global lower bound
must differ by at most
\(10^{-7}\max\{1,|\widehat v(q^{\mathrm{SAA}})|\}\), up to numerical
roundoff. Independent recourse evaluation checks the reported cost.

The reported construction times were measured on an Intel Core
i7-9750H at 2.60 GHz, running Ubuntu 26.04 LTS under WSL2, with
Python 3.12.13, NumPy 2.3.5, SciPy 1.17.0, and HiGHS 1.8.0.
These measurements use CPU solvers, one numerical-library thread,
and sequential execution after the parameter sweep.

\subsection{Uniform-metric validation}
\label{app:numerics-uniform}

We use \(d_{ij}=1\) for \(i\ne j\) to check the implementation against
Proposition~\ref{prop:ih-uniform-exactness}.
Table~\ref{tab:app-uniform} aggregates 800 completed instances at each
transportation scale over network sizes, demand laws, CVs, and
critical fractiles.

\begin{table}[htbp]
\centering
\small
\caption{Uniform-metric validation, \(N=2000\). Inventory ratios are
averaged across instances; residual columns report maxima.}
\label{tab:app-uniform}
\resizebox{\textwidth}{!}{%
\begin{tabular}{rrrrrr}
\toprule
\(\gamma\) & \(N\) & Mean \(\Delta(\mathrm{ALG})\) (\%)
& \(Q^{\mathrm{ALG}}/Q^{\mathrm{SAA}}\)
& Equation residual & Raw CDF residual \\
\midrule
0.5 & 2000 & 0.0251 & 1.0000 & \(<10^{-9}\) & 0.0002249 \\
1.0 & 2000 & 0.0367 & 1.0000 & \(<10^{-9}\) & 0.0001988 \\
2.0 & 2000 & 0.0512 & 1.0001 & \(<10^{-9}\) & 0.0001499 \\
3.0 & 2000 & 0.0634 & 1.0000 & \(<10^{-9}\) & 0.0001000 \\
4.5 & 2000 & 0.0795 & 1.0003 & \(<10^{-9}\) & 0.0000250 \\
\bottomrule
\end{tabular}%
}
\end{table}

Mean sampled gaps remain below \(0.08\%\), and mean total-inventory
ratios are close to one. The population exactness result uses exact
distributional CDFs. Here, analytic marginal CDFs are combined with
interpolated empirical pooling CDFs, whereas SAA optimizes the sampled
objective directly. A nonzero sampled gap therefore does not
contradict the uniform-metric characterization.

\subsection{Cost decomposition}
\label{app:numerics-costs}

Using \eqref{eq:model-decomp}, the reported components are
\[
\begin{aligned}
\text{Holding}
&=\frac hN\sum_{s=1}^{N}(Q-D^{(s)})_+,\\
\text{Backlog}
&=\frac bN\sum_{s=1}^{N}(D^{(s)}-Q)_+,\\
\text{Transportation}
&=\frac1N\sum_{s=1}^{N}T(q,\xi^{(s)}).
\end{aligned}
\]
Their sum equals \(\widehat v(q)\). Table~\ref{tab:app-cost-components}
uses the same paired instances as Table~\ref{tab:numerics-main}.

\begin{table}[htbp]
\centering
\small
\caption{Mean cost components on completed paired instances,
\(N=2000\). Displayed components may differ from totals by rounding.}
\label{tab:app-cost-components}
\resizebox{\textwidth}{!}{%
\begin{tabular}{llrrrrrr}
\toprule
Metric & Policy & \(N\) & Paired cases
& Holding & Backlog & Transportation & Total cost \\
\midrule
Star & SAA & 2000 & 3998 & 3.171 & 0.237 & 3.151 & 6.559 \\
Star & IND & 2000 & 3998 & 6.663 & 0.037 & 2.149 & 8.848 \\
Star & ALG & 2000 & 3998 & 3.122 & 0.322 & 3.800 & 7.244 \\
\midrule
Line & SAA & 2000 & 3486 & 1.554 & 0.542 & 0.521 & 2.617 \\
Line & IND & 2000 & 3486 & 6.737 & 0.034 & 0.126 & 6.897 \\
Line & ALG & 2000 & 3486 & 2.113 & 0.269 & 0.475 & 2.857 \\
\midrule
Clustered & SAA & 2000 & 3573 & 1.671 & 0.421 & 0.978 & 3.071 \\
Clustered & IND & 2000 & 3573 & 6.913 & 0.026 & 0.278 & 7.217 \\
Clustered & ALG & 2000 & 3573 & 1.816 & 0.324 & 0.994 & 3.134 \\
\midrule
Geometric & SAA & 2000 & 7369 & 1.407 & 0.461 & 0.800 & 2.667 \\
Geometric & IND & 2000 & 7369 & 4.839 & 0.051 & 0.283 & 5.173 \\
Geometric & ALG & 2000 & 7369 & 1.531 & 0.380 & 0.913 & 2.823 \\
\midrule
Random graph & SAA & 2000 & 11816 & 1.953 & 0.365 & 2.543 & 4.861 \\
Random graph & IND & 2000 & 11816 & 6.178 & 0.037 & 1.114 & 7.328 \\
Random graph & ALG & 2000 & 11816 & 2.002 & 0.337 & 2.587 & 4.925 \\
\bottomrule
\end{tabular}%
}
\end{table}

Relative to IND, ALG lowers holding costs while increasing backlog
and transportation costs in every metric class. The holding-cost
reduction dominates these increases, producing the average savings
reported in the main text. Compared with SAA, the star instances show
the largest difference in mean transportation cost, consistent with
the different placements illustrated in
Figure~\ref{fig:numerics-allocation}. Holding and backlog depend on
total inventory, whereas transportation also depends on its placement.

\subsection{Critical fractile and demand sensitivity}
\label{app:numerics-demand}

Tables~\ref{tab:app-fractile-cv} and~\ref{tab:app-demand-law} separate
ALG performance by critical fractile, CV, and demand law, averaging
over the remaining parameters of the main grid.

\begin{table}[htbp]
\centering
\small
\caption{Mean ALG sampled gap (\%) by critical fractile and CV,
\(N=2000\).}
\label{tab:app-fractile-cv}
\resizebox{\textwidth}{!}{%
\begin{tabular}{rrrrrrrr}
\toprule
\(b/(b+h)\) & CV & \(N\)
& Star & Line & Clustered & Geometric & Random graph \\
\midrule
0.80 & 0.25 & 2000 & 6.36 & 5.99 & 5.61 & 6.77 & 4.29 \\
0.80 & 0.50 & 2000 & 7.15 & 4.58 & 1.93 & 2.62 & 1.57 \\
0.85 & 0.25 & 2000 & 7.69 & 9.42 & 3.40 & 8.40 & 4.60 \\
0.85 & 0.50 & 2000 & 8.48 & 6.55 & 1.61 & 3.55 & 1.29 \\
0.90 & 0.25 & 2000 & 8.66 & 13.24 & 2.40 & 9.43 & 3.12 \\
0.90 & 0.50 & 2000 & 9.46 & 9.97 & 1.49 & 4.21 & 1.22 \\
0.95 & 0.25 & 2000 & 8.75 & 18.51 & 2.06 & 9.36 & 1.13 \\
0.95 & 0.50 & 2000 & 9.85 & 13.85 & 1.76 & 4.91 & 0.79 \\
\bottomrule
\end{tabular}%
}
\end{table}

Larger critical fractiles increase the reported mean gaps on stars
and lines. On lines with \(\mathrm{CV}=0.25\), the gap rises from
\(5.99\%\) at \(b/(b+h)=0.80\) to \(18.51\%\) at \(0.95\).
Higher variability does not uniformly increase the gap: at each
tested fractile, \(\mathrm{CV}=0.50\) produces smaller mean gaps
than \(0.25\) on lines, clustered instances, geometric instances,
and random graphs, while the reverse holds on stars.

\begin{table}[htbp]
\centering
\small
\caption{Mean ALG sampled gap (\%) by demand law, \(N=2000\).}
\label{tab:app-demand-law}
\resizebox{\textwidth}{!}{%
\begin{tabular}{lrrrrrr}
\toprule
Demand law & \(N\)
& Star & Line & Clustered & Geometric & Random graph \\
\midrule
Truncated Gaussian & 2000 & 5.85 & 11.92 & 1.79 & 7.65 & 2.86 \\
Log-normal & 2000 & 8.23 & 10.68 & 1.69 & 6.91 & 2.43 \\
Gamma & 2000 & 7.35 & 11.46 & 2.03 & 7.26 & 2.36 \\
Pareto & 2000 & 11.75 & 5.60 & 5.13 & 3.93 & 1.22 \\
\bottomrule
\end{tabular}%
}
\end{table}

The ordering of demand laws differs across metrics. Pareto demand
has the largest mean gap on stars and clustered instances, but the
smallest on lines, geometric instances, and random graphs. Thus
matching the mean and CV does not remove the influence of distributional
shape, and no demand law gives the largest reported gap across all
metric classes.

\subsection{Metric sensitivity}
\label{app:numerics-metrics}

For clustered instances, we fix \(n=30\), divide the nodes into two
equal groups, and compare within-group base distances
\(0.1\), \(0.25\), and \(0.5\). Between-group distance remains \(1\).
The baseline setting is included once, and all policy statistics in
each row use the same paired instances.

\begin{table}[htbp]
\centering
\small
\caption{Cluster-separation sensitivity at \(n=30\), \(N=2000\).
Other parameters follow the common grid.}
\label{tab:app-cluster}
\resizebox{\textwidth}{!}{%
\begin{tabular}{rrrrrrr}
\toprule
Within-group distance & \(N\) & Paired cases
& Mean \(\Delta(\mathrm{IND})\) (\%)
& Mean \(\Delta(\mathrm{ALG})\) (\%)
& P90 \(\Delta(\mathrm{ALG})\) (\%)
& ALG imp.\ vs.\ IND (\%) \\
\midrule
0.10 & 2000 & 736 & 145.17 & 2.62 & 4.20 & 54.12 \\
0.25 & 2000 & 783 & 102.39 & 0.80 & 0.90 & 45.24 \\
0.50 & 2000 & 800 & 66.78 & 0.12 & 0.26 & 33.89 \\
\bottomrule
\end{tabular}%
}
\end{table}

As within-group costs approach between-group costs, the metric
becomes closer to uniform. Over these settings, the mean ALG gap
declines from \(2.62\%\) to \(0.12\%\). IND also becomes closer
to SAA, so the mean improvement of ALG over IND decreases from
\(54.12\%\) to \(33.89\%\).

For random graphs, we fix \(n=30\) and vary
\(p\in\{0.2,0.4,0.6\}\). Graphs are conditioned on connectivity
and use unit-edge shortest-path distances, normalized so that the
maximum transportation cost is \(\gamma\).

\begin{table}[htbp]
\centering
\small
\caption{Random-graph density sensitivity at \(n=30\), \(N=2000\).
Other parameters follow the common grid.}
\label{tab:app-density}
\resizebox{\textwidth}{!}{%
\begin{tabular}{rrrrrrr}
\toprule
\(p\) & \(N\) & Paired cases
& Mean \(\Delta(\mathrm{IND})\) (\%)
& Mean \(\Delta(\mathrm{ALG})\) (\%)
& P90 \(\Delta(\mathrm{ALG})\) (\%)
& ALG imp.\ vs.\ IND (\%) \\
\midrule
0.20 & 2000 & 2015 & 95.06 & 3.09 & 7.94 & 41.96 \\
0.40 & 2000 & 2206 & 79.79 & 2.84 & 6.48 & 37.24 \\
0.60 & 2000 & 2294 & 54.14 & 1.98 & 2.57 & 28.95 \\
\bottomrule
\end{tabular}%
}
\end{table}

The reported mean ALG gap decreases from \(3.09\%\) to \(1.98\%\)
as \(p\) increases from \(0.2\) to \(0.6\). The improvement over
IND also decreases, from \(41.96\%\) to \(28.95\%\).
Because each graph is normalized by its realized diameter, this
comparison changes both path structure and relative transportation
costs while keeping the maximum cost fixed.

\subsection{Computational performance}
\label{app:numerics-runtime}

Table~\ref{tab:app-runtime} reports cold-start construction times
for the random-graph groups at \(n=40\) and \(n=50\), where all
eight settings yield completed solutions for every policy. Each
group comprises the four demand laws and two CV values, with
\(p=0.4\), \(b/(b+h)=0.90\), \(\gamma=2\), and \(N=2000\).
ALG time includes recourse calibration, weighted-CDF construction,
and solution of the inventory equations. SAA time includes model
construction and optimization. Common cost evaluation is excluded,
and every timing uses a fresh policy construction.

\begin{table}[htbp]
\centering
\small
\caption{Mean policy-construction time in seconds for completed
random-graph groups, \(N=2000\).}
\label{tab:app-runtime}
\begin{tabular}{rrrrrr}
\toprule
\(n\) & \(N\) & Paired settings & IND & ALG & SAA \\
\midrule
40 & 2000 & 8 & 0.00023 & 0.571 & 3.693 \\
50 & 2000 & 8 & 0.00022 & 0.927 & 6.507 \\
\bottomrule
\end{tabular}
\end{table}

The ratio of mean SAA to ALG construction time is approximately
\(6.5\) at \(n=40\) and \(7.0\) at \(n=50\). These comparisons
hold the sample size fixed at \(N=2000\).

\end{document}